\documentclass[12pt]{article}

\usepackage{amsmath, epsfig, cite}
\usepackage{amssymb}
\usepackage{amsfonts}
\usepackage{latexsym}
\usepackage{float}
\usepackage{booktabs}
\usepackage{caption}

\newtheorem{thm}{Theorem}[section]

\newtheorem{cor}[thm]{Corollary}

\newtheorem{lem}[thm]{Lemma}

\numberwithin{equation}{section}

\UseRawInputEncoding
\begin{document}


\begin{center}
{\Large\bf On the positive coefficients of two families of $q$-series}
\end{center}

\vskip 2mm \centerline{Ji-Cai Liu and Kong-Lian Liao}
\begin{center}
{\footnotesize Department of Mathematics, Wenzhou University, Wenzhou 325035, PR China\\
{\tt jcliu2016@gmail.com, 23450829013@stu.wzu.edu.cn} }
\end{center}


\vskip 0.7cm \noindent{\bf Abstract.}
Let $S$ be a finite set of pairwise coprime positive integers and $Ax^2+Bx$ be an integer valued polynomial with $A> B\ge 0$. For integers $k\ge 1$ and $n\ge 0$, the coefficients $\gamma_{S,A,B}^k (n)$ are defined as
\begin{align*}
\prod_{s\in S}\frac{1}{1-q^s}\sum_{j\not\in [-k,k-1]} (-1)^{j+k}q^{Aj^2+Bj}=\sum_{n= 0}^{\infty}\gamma_{S,A,B}^k (n)q^n.
\end{align*}
In this paper, we investigate the positivity of $\gamma_{S,A,B}^k (n)$ for $|S|=4,5$.

\vskip 3mm \noindent {\it Keywords}: positivity; partitions; pentagonal number theorem

\vskip 2mm
\noindent{\it MR Subject Classifications}: 11S05, 05A17, 05A20
	
\section{Introduction}
Euler's pentagonal number theorem \cite[Corollary 1.7]{andrews-b-1998} plays an important role in the theory of integer partitions, which is stated as follows:
\begin{align}
(q;q)_{\infty}=\sum_{j=-\infty}^{\infty}(-1)^jq^{j(3j+1)/2}.\label{a-1}
\end{align}
Here and throughout this paper, the $q$-shifted factorial is defined by $(a;q)_0=1$,
$(a;q)_n=(1-a)(1-aq)\cdots (1-aq^{n-1})$ for $n\ge 1$, $(a;q)_{\infty}=\prod_{k=0}^{\infty}(1-aq^k)$ and
$(a_1,a_2,\cdots,a_m;q)_{\infty}=(a_1;q)_{\infty}(a_2;q)_{\infty}\cdots (a_m;q)_{\infty}$.
The $q$-binomial coefficient is defined as
\begin{align*}
{n\brack k}={n\brack k}_{q}
=\begin{cases}
\displaystyle\frac{(q;q)_n}{(q;q)_k(q;q)_{n-k}} &\text{if $0\le k\le n$},\\[10pt]
0 &\text{otherwise.}
\end{cases}
\end{align*}

Andrews and Merca \cite{am-jcta-2012} proved that for $k\ge 1$,
\begin{align}
(-1)^{k-1}\sum_{j=0}^{k-1}(-1)^j\left(p(n-j(3j+1)/2)-p(n-(j+1)(3j+2)/2)\right)\ge 0,\label{a-2}
\end{align}
where $p(n)$ denotes the number of all partitions of $n$.
Note that \eqref{a-2} is equivalent to \begin{align}
\frac{(-1)^{k-1}}{(q;q)_{\infty}}\sum_{j=0}^{k-1} (-1)^jq^{j(3j+1)/2}(1-q^{2j+1})+(-1)^k
\in \mathbb{N}[[q]].\label{a-3}
\end{align}
By \eqref{a-1}, we can rewrite \eqref{a-3} as
\begin{align}
\frac{1}{(q;q)_{\infty}}\sum_{j\not \in [-k,k-1]}(-1)^{j+k}q^{j(3j+1)/2}
\in \mathbb{N}[[q]].\label{a-4}
\end{align}
Here and throughout this paper, we use the notation: for integers $a$ and $b$ with $a<b$,
\begin{align*}
\sum_{j\not\in[a,b]} A_j=\sum_{j=-\infty}^{a-1}A_j+\sum_{j=b+1}^{\infty}A_j.
\end{align*}

Let $S$ be a finite set of pairwise coprime positive integers and $Ax^2+Bx$ be an integer valued polynomial with $A> B\ge 0$. For integers $k\ge 1$ and $n\ge 0$, the coefficients $\gamma_{S,A,B}^k (n)$ are defined as
\begin{align*}
\prod_{s\in S}\frac{1}{1-q^s}\sum_{j\not\in [-k,k-1]} (-1)^{j+k}q^{Aj^2+Bj}=\sum_{n= 0}^{\infty}\gamma_{S,A,B}^k (n)q^n.
\end{align*}

Yao \cite{yao-jcta-2024} strengthened \eqref{a-4} as follows:
\begin{align*}
\gamma_{\{1,2,3\},3/2,1/2}^k (n)\ge 0\quad\text{for $k\ge 1$ and $n\ge 0$.}
\end{align*}
Subsequently, Zhou \cite{zhou-jcta-2024} studied the positivity of $\gamma_{\{a,b,c\},3/2,1/2}^k (n)$
for pairwise coprime positive integers $a,b,c$ and found the following triples $(a,b,c)$ such that $\gamma_{\{a,b,c\},3/2,1/2}^k (n)\ge 0$ for $k\ge 1$ and $n\ge 0$:
\begin{align*}
(a,b,c)\in\{(1,2,3),(1,2,5),(1,2,7),(1,3,4),(1,3,5)\}.
\end{align*}

The first author \cite{liu-a-2024} investigated the positivity of $\gamma_{\{a,b,c\},A,B}^k (n)$ for any
pairwise coprime positive integers $a,b,c$ and integer valued polynomial $Ax^2+Bx$ with $A> B\ge 0$. The first author \cite{liu-a-2024} determined the two values $K_{\{a,b,c\},A,B}$ and $N_{\{a,b,c\},A,B}^k$ such that
\begin{align}
&\gamma_{\{a,b,c\},A,B}^k (n)\ge 0\quad\text{for $n\ge 0$ and $k\ge K_{\{a,b,c\},A,B}$},\label{a-5}\\[5pt]
&\gamma_{\{a,b,c\},A,B}^k (n)\ge 0\quad \text{for $n\ge N_{\{a,b,c\},A,B}^k$}.\notag
\end{align}

For $|S|\ge 4$, let $S'$ denote the set of the smallest three elements in $S$. Note that
\begin{align}
\sum_{n= 0}^{\infty}\gamma_{S,A,B}^k (n)q^n&=\prod_{s\in S\setminus S'}\frac{1}{1-q^s}\prod_{s\in S'}\frac{1}{1-q^s}\sum_{j\not\in [-k,k-1]} (-1)^{j+k}q^{Aj^2+Bj}\notag\\[5pt]
&=\prod_{s\in S\setminus S'}\frac{1}{1-q^s}\sum_{n= 0}^{\infty}\gamma_{S',A,B}^k (n)q^n.\label{a-6}
\end{align}
By \eqref{a-5} and \eqref{a-6}, for $k\ge K_{S',A,B}$ and $n\ge 0$, we have $\gamma_{S',A,B}^k (n)\ge 0$
and so $\gamma_{S,A,B}^k (n)\ge 0$.

Sometimes we need a larger set $S$ to ensure that $\gamma_{S,A,B}^k (n)\ge 0$ for all integers $k\ge 1$ and $n\ge 0$.
The starting point of the paper is to determine the value $N_{S,A,B}^k$ for $|S|=4,5$ such that $\gamma_{S,A,B}^k (n)\ge 0$ for $n\ge N_{S,A,B}^k$.
Once the value $N_{S,A,B}^k$ is determined for $|S|=4,5$, we only need verify a finite number of $\gamma_{S,A,B}^k (n)$ for $n<N_{S,A,B}^k$ with $k<K_{S',A,B}$. As an application, we reprove two positivity conjectures due to Merca \cite{merca-em-2021}.

The rest of the paper is organized as follows. We state the main theorems and some corollaries
in Section 2. Section 3 is devoted to some preliminary results. The proofs of the main theorems are presented in Sections 4 and 5.

\section{Main results}
For a real polynomial $G(x)$, let $\mathcal{R}_xG(x)$ be the set of all real roots of $G(x)$.
For real polynomials $G_1(x),G_2(x),\cdots,G_m(x)$, let
\begin{align*}
\mathcal{T}_x\{G_1(x),G_2(x),\cdots,G_m(x)\}=\max\left(\{1\}\cup\mathcal{R}_xG_1(x)\cup \mathcal{R}_xG_2(x)\cup\cdots \cup \mathcal{R}_xG_m(x)\right).
\end{align*}
For a real number $x$, let $\lceil x\rceil$ denote the smallest integer greater than or equal to $x$.

Now we are ready to state the main theorems.
\begin{thm}\label{t-1}
Let $Y=\{y_1,y_2,y_3,y_4\}$ in which $y_1,y_2,y_3,y_4$ are pairwise coprime positive integers.
Then $\gamma_{Y,A,B}^k (n)\ge 0$ for $n\ge N_{Y,A,B}^k=A\left(k+2L_{Y,A,B}^k\right)^2+B\left(k+2L_{Y,A,B}^k\right)$
with
\begin{align*}
L_{Y,A,B}^k=\lceil \mathcal{T}_l\{H_1,H_2,\cdots,H_8\}\rceil,
\end{align*}
where $H_1,H_2,\cdots,H_8$ are listed in Appendix A.
\end{thm}

\begin{thm}\label{t-2}
Let $Z=\{z_1,z_2,z_3,z_4,z_5\}$ in which $z_1,z_2,z_3,z_4,z_5$ are pairwise coprime positive integers.
Then $\gamma_{Z,A,B}^k (n)\ge 0$ for $n\ge N_{Z,A,B}^k=A\left(k+2L_{Z,A,B}^k\right)^2+B\left(k+2L_{Z,A,B}^k\right)$ with
\begin{align*}
L_{Z,A,B}^k=\lceil \mathcal{T}_l\{G_1,G_2,\cdots,G_{12}\}\rceil,
\end{align*}
where $G_1,G_2,\cdots,G_{12}$ are listed in Appendix B.
\end{thm}

Merca \cite{merca-em-2021} proposed the following two positivity conjectures (equivalent forms):
\begin{align}
&\frac{1}{(q,q^4;q^5)_{\infty}}\sum_{j\not\in [-k,k-1]} (-1)^{j+k}q^{Aj^2+Bj} \in \mathbb{N}[[q]],\label{m-1}\\[5pt]
&\frac{1}{(q^2,q^3;q^5)_{\infty}}\sum_{j\not\in [-k,k-1]} (-1)^{j+k}q^{Aj^2+Bj} \in \mathbb{N}[[q]],\label{m-2}
\end{align}
for $(A,B)\in\{(5/2,3/2),(5/2,1/2),(7/2,5/2),(7/2,3/2)\}$. Chen and Yao \cite{cy-jcta-2024} proved
\eqref{m-1} and \eqref{m-2} by establishing the following two corollaries:

\begin{cor}
Let $Z=\{1,4,9,11,19\}$ and
\begin{align*}
(A,B)\in\{(5/2,3/2),(5/2,1/2),(7/2,5/2),(7/2,3/2)\}.
\end{align*}
Then $\gamma_{Z,A,B}^k (n)\ge 0$ for all integers $k\ge 1$ and $n\ge 0$.
\end{cor}

\begin{cor}
Let $Z=\{2,3,7,13,17\}$ and
\begin{align*}
(A,B)\in\{(5/2,3/2),(5/2,1/2),(7/2,5/2),(7/2,3/2)\}.
\end{align*}
Then $\gamma_{Z,A,B}^k (n)\ge 0$ for all integers $k\ge 1$ and $n\ge 0$.
\end{cor}

By Theorem \ref{t-2} and the result due to the first author \cite{liu-a-2024}, we can compute the values
$K_{Z',A,B}$ and $N_{Z,A,B}^k$. To prove the above two corollaries, it suffices to verify a finite number of $\gamma_{Z,A,B}^k (n)$ for $n<N_{Z,A,B}^k$ with $k<K_{Z',A,B}$.

\begin{table}[H]
\caption*{\text{The case $Z=\{1,4,9,11,19\}$}}
\centering
\begin{tabular}{ccc}
\toprule
$(A,B)$&$K_{Z',A,B}$&$\left\{N_{Z,A,B}^k\right\}_{k=1}^{K_{Z',A,B}-1}$\\
\midrule
$(5/2,3/2)$&$5$&$216,99,70,99$\\[5pt]
$(5/2,1/2)$&$4$&$207,93,65$\\[5pt]
$(7/2,5/2)$&$4$&$189,141,100$\\[5pt]
$(7/2,3/2)$&$3$&$182,62$\\
\bottomrule
\end{tabular}
\end{table}

\begin{table}[H]
\caption*{\text{The case $Z=\{2,3,7,13,17\}$}}
\centering
\begin{tabular}{ccc}
\toprule
$(A,B)$&$K_{Z',A,B}$&$\left\{N_{Z,A,B}^k\right\}_{k=1}^{K_{Z',A,B}-1}$\\
\midrule
$(5/2,3/2)$&$5$&$216,99,70,99$\\[5pt]
$(5/2,1/2)$&$4$&$207,93,65$\\[5pt]
$(7/2,5/2)$&$4$&$189,141,100$\\[5pt]
$(7/2,3/2)$&$3$&$182,62$\\
\bottomrule
\end{tabular}
\end{table}

In fact, we can also reprove \eqref{m-1} and \eqref{m-2} by establishing the following two corollaries.
\begin{cor}
Let $Z=\{1,4,9,11,29\}$ and
\begin{align*}
(A,B)\in\{(5/2,3/2),(5/2,1/2),(7/2,5/2),(7/2,3/2)\}.
\end{align*}
Then $\gamma_{Z,A,B}^k (n)\ge 0$ for all integers $k\ge 1$ and $n\ge 0$.
\end{cor}

\begin{cor}
Let $Z=\{2,3,7,13,23\}$ and
\begin{align*}
(A,B)\in\{(5/2,3/2),(5/2,1/2),(7/2,5/2),(7/2,3/2)\}.
\end{align*}
Then $\gamma_{Z,A,B}^k (n)\ge 0$ for all integers $k\ge 1$ and $n\ge 0$.
\end{cor}

\begin{table}[H]
\caption*{\text{The case $Z=\{1,4,9,11,29\}$}}
\centering
\begin{tabular}{ccc}
\toprule
$(A,B)$&$K_{Z',A,B}$&$\left\{N_{Z,A,B}^k\right\}_{k=1}^{K_{Z',A,B}-1}$\\
\midrule
$(5/2,3/2)$&$5$&$442,172,133,99$\\[5pt]
$(5/2,1/2)$&$4$&$429,164,65$\\[5pt]
$(7/2,5/2)$&$4$&$451,141,100$\\[5pt]
$(7/2,3/2)$&$3$&$297,135$\\
\bottomrule
\end{tabular}
\end{table}

\begin{table}[H]
\caption*{\text{The case $Z=\{2,3,7,13,23\}$}}
\centering
\begin{tabular}{ccc}
\toprule
$(A,B)$&$K_{Z',A,B}$&$\left\{N_{Z,A,B}^k\right\}_{k=1}^{K_{Z',A,B}-1}$\\
\midrule
$(5/2,3/2)$&$5$&$319,172,133,99$\\[5pt]
$(5/2,1/2)$&$4$&$308,93,65$\\[5pt]
$(7/2,5/2)$&$4$&$306,141,100$\\[5pt]
$(7/2,3/2)$&$3$&$182,62$\\
\bottomrule
\end{tabular}
\end{table}

\section{Preliminaries}
\begin{lem} (See \cite[(3.6)]{liu-a-2024}.) \label{lem-1}
For pairwise coprime positive integers $y_1,y_2,y_3$ and $y_4$, we have
\begin{align*}
\frac{1}{(1-q^{y_1})(1-q^{y_2})(1-q^{y_3})(1-q^{y_4})}=\sum_{n\ge 0}(an^3+bn^2+cn)q^n+
\sum_{n\ge 0}s(n)q^n,
\end{align*}
where $\{s(n)\}_{n\ge 0}$ has a period of $y_1y_2y_3y_4$ and
\begin{align*}
a&=\frac{1}{6y_1y_2y_3y_4},\\[5pt]
b&=\frac{y_1+y_2+y_3+y_4}{4y_1y_2y_3y_4},\\[5pt]
c&=\frac{y_1^2+y_2^2+y_3^2+y_4^2+3(y_1y_2+y_1y_3+y_1y_4+y_2y_3+y_2y_4+y_3y_4)}{12y_1y_2y_3y_4}.
\end{align*}
\end{lem}
Since $\{s(n)\}_{n\ge 0}$ is a periodic sequence, there exists a smallest bound $D_{y_1,y_2,y_3,y_4}$ such that $|s(n)|\le D_{y_1,y_2,y_3,y_4}$ for all integers $n\ge 0$.

\begin{lem} (See \cite[(3.7)]{liu-a-2024}.)\label{lem-2}
For pairwise coprime positive integers $z_1,z_2,z_3,z_4$ and $z_5$, we have
\begin{align*}
\frac{1}{(1-q^{z_1})(1-q^{z_2})(1-q^{z_3})(1-q^{z_4})(1-q^{z_5})}=\sum_{n\ge 0}(an^4+bn^3+cn^2+dn)q^n+
\sum_{n\ge 0}t(n)q^n,
\end{align*}
where $\{t(n)\}_{n\ge 0}$ has a period of $z_1z_2z_3z_4z_5$ and
\begin{align*}
a&=\frac{1}{24z_1z_2z_3z_4z_5},\\[5pt]
b&=\frac{z_1+z_2+z_3+z_4+z_5}{12z_1z_2z_3z_4z_5},\\[5pt]
c&=\frac{z_1^2+z_2^2+z_3^2+z_4^2+z_5^2}{24z_1z_2z_3z_4z_5}\\[5pt]
&+\frac{z_1z_2+z_1z_3+z_1z_4+z_1z_5+z_2z_3+z_2z_4+z_2z_5+z_3z_4+z_3z_5+z_4z_5}{8z_1z_2z_3z_4z_5},\\[5pt]
d&=\frac{(z_1+z_2+z_3+z_4+z_5)}{24z_1z_2z_3z_4z_5}\\[5pt]
&\times (z_1z_2+z_1z_3+z_1z_4+z_1z_5+z_2z_3+z_2z_4+z_2z_5+z_3z_4+z_3z_5+z_4z_5).
\end{align*}
\end{lem}
Since $\{t(n)\}_{n\ge 0}$ is a periodic sequence, there exists a smallest bound $D_{z_1,z_2,z_3,z_4,z_5}$ such that $|t(n)|\le D_{z_1,z_2,z_3,z_4,z_5}$ for all integers $n\ge 0$.

We remark that P\'olya and Szeg\H{o} \cite[Problem 27.1, page 5]{ps-b-1998} gave a general result of Lemmas \ref{lem-1} and \ref{lem-2} without explicit polynomial formulas.

We also require the following trivial result.
\begin{lem}\label{lem-3}
Let $G_1(x),G_2(x),\cdots,G_m(x)$ be real polynomials with positive leading coefficients. If
$x\ge \mathcal{T}_x\{G_1(x),G_2(x),\cdots,G_m(x)\}$, then $G_i(x)\ge 0$ for $i=1,2,\cdots,m$.
\end{lem}

\section{Proof of Theorems \ref{t-1}}
Let $a,b,c$ be the quantities defined in Lemma \ref{lem-1} and $e=D_{y_1,y_2,y_3,y_4}$.
Let $F(x)=ax^3+bx^2+cx$, $f(j)=Aj^2+Bj$ and $g(j)=Aj^2-Bj$.
It is trivial to check that for all integers $k\ge 1$ and $j\ge 0$,
\begin{align*}
f(k+2j)<g(k+2j+1)\le f(k+2j+1)<g(k+2j+2)\le f(k+2j+2).
\end{align*}
For any integer $n\ge f(k)$, there exists a unique integer $l\ge 0$ such that
\begin{align*}
f(k+2l)\le n< f(k+2l+2).
\end{align*}
For a real polynomial $G(x)$, let $LC_x(G(x))$ denote the leading coefficient of $G(x)$.

Next, we shall distinguish four cases to prove Theorem \ref{t-1}.

{\noindent\bf Case 1} $f(k+2l)\le n<g(k+2l+1)$ with $l\ge 1$. We have
\begin{align*}
\gamma_{Y,A,B}^k (n)&\ge \sum_{j=0}^{l-1}F(n-f(k+2j))-\sum_{j=0}^{l-1}F(n-g(k+2j+1))\\[5pt]
&-\sum_{j=0}^{l-1}F(n-f(k+2j+1))+\sum_{j=0}^{l-1}F(n-g(k+2j+2))\\[5pt]
&+F(n-f(k+2l))-(4l+1)e\\[5pt]
&=P_1(n).
\end{align*}
With the help of Maple, we find that $P_1(n)$ is a cubic polynomial in variable $n$ with
$LC_n(P_1(n))=a>0$ and
\begin{align*}
P_2(n)=\frac{\mathrm{d}}{\mathrm{d}n}P_1(n)=3an^2+Un+V,
\end{align*}
where $3a>0$, $U$ and $V$ are independent of $n$. It follows that
\begin{align*}
P_2(n)\ge P_2\left(-\frac{U}{6a}\right)=H_1.
\end{align*}
Note that $LC_l(H_1)=12a(A-B)(2Ak-A-B)>0$. By Lemma \ref{lem-3}, we have $P_2(n)\ge H_1\ge 0$ for $l\ge \mathcal{T}_l\{H_1\}$. For $l\ge \mathcal{T}_l\{H_1\}$, $P_1(n)$ is a weakly increasing function on the interval $[f(k+2l), g(k+2l+1))$, and so
\begin{align*}
\gamma_{Y,A,B}^k (n)\ge P_1(n)\ge P_1(f(k+2l))=H_2.
\end{align*}
Noting that $LC_l(H_2)=48A^2ak(A-B)>0$, by Lemma \ref{lem-3} we have
$\gamma_{Y,A,B}^k (n)\ge H_2\ge 0$ for $l\ge \mathcal{T}_l\{H_1,H_2\}$.

Then $\gamma_{Y,A,B}^k (n)\ge 0$ for $f(k+2l)\le n<g(k+2l+1)$ with $l\ge \mathcal{T}_l\{H_1,H_2\}$.
\\[10pt]
{\noindent\bf Case 2} $g(k+2l+1)\le n<f(k+2l+1)$ with $l\ge 1$. We have
\begin{align*}
\gamma_{Y,A,B}^k (n)&\ge \sum_{j=0}^{l-1}F(n-f(k+2j))-\sum_{j=0}^{l-1}F(n-g(k+2j+1))\\[5pt]
&-\sum_{j=0}^{l-1}F(n-f(k+2j+1))+\sum_{j=0}^{l-1}F(n-g(k+2j+2))\\[5pt]
&+F(n-f(k+2l))-F(n-g(k+2l+1))-(4l+2)e\\[5pt]
&=P_3(n).
\end{align*}
Via Maple, we find that
\begin{align*}
P_3(n)=3 a (2 k + 2 l + 1) (A - B)n^2+Un+V,
\end{align*}
where $3 a (2 k + 2 l + 1) (A - B)>0$, $U$ and $V$ are independent of $n$. It follows that
\begin{align*}
\gamma_{Y,A,B}^k (n)\ge P_3(n)\ge P_3\left(-\frac{U}{6a (2 k + 2 l + 1) (A - B)}\right)=H_3.
\end{align*}
Noting that $LC_l(H_3)=48 A^2 a k (A-B)>0$, by Lemma \ref{lem-3} we have
$\gamma_{Y,A,B}^k (n)\ge H_3\ge 0$ for $l\ge \mathcal{T}_l \{H_3\}$.

Then $\gamma_{Y,A,B}^k (n)\ge 0$ for $g(k+2l+1)\le n<f(k+2l+1)$ with $l\ge \mathcal{T}_l \{H_3\}$.
\\[10pt]
{\noindent\bf Case 3} $f(k+2l+1)\le n<g(k+2l+2)$ with $l\ge 1$. We have
\begin{align*}
\gamma_{Y,A,B}^k (n)&\ge \sum_{j=0}^{l-1}F(n-f(k+2j))-\sum_{j=0}^{l-1}F(n-g(k+2j+1))\\[5pt]
&-\sum_{j=0}^{l-1}F(n-f(k+2j+1))+\sum_{j=0}^{l-1}F(n-g(k+2j+2))\\[5pt]
&+F(n-f(k+2l))-F(n-g(k+2l+1))-F(n-f(k+2l+1))-(4l+3)e\\[5pt]
&=P_4(n).
\end{align*}
Via Maple, we find that $P_4(n)$ is a cubic polynomial in variable $n$ with
$LC_n(P_4(n))=-a<0$ and
\begin{align*}
P_5(n)=\frac{\mathrm{d}}{\mathrm{d}n}P_4(n)=-3an^2+Un+V,
\end{align*}
where $-3a<0$, $U$ and $V$ are independent of $n$. It follows that
\begin{align*}
P_5(n)\ge \min\{P_5(f(k+2l+1)),P_5(g(k+2l+2))\}.
\end{align*}
Let $H_4=P_5(f(k+2l+1))$ and $H_5=P_5(g(k+2l+2))$. Since $LC_l(H_4)=LC_l(H_5)=24a(A-B)(Ak+B)>0$, by Lemma \ref{lem-3} we have $P_5(n)\ge \min\{H_4,H_5\}\ge 0$ for $l\ge \mathcal{T}\{H_4,H_5\}$. It follows that
for $l\ge \mathcal{T}\{H_4,H_5\}$, $P_4(n)$ is a weakly increasing function on the interval $[f(k+2l+1), g(k+2l+2))$, and so
\begin{align*}
\gamma_{Y,A,B}^k (n)\ge P_4(n)\ge P_4(f(k+2l+1))=H_6.
\end{align*}
Noting that $LC_l(H_6)=48A^2ak(A-B)>0$, by Lemma \ref{lem-3} we have $\gamma_{Y,A,B}^k (n)\ge H_6\ge 0$ for $l\ge \mathcal{T}\{H_4,H_5,H_6\}$.

Then $\gamma_{Y,A,B}^k (n)\ge 0$ for $f(k+2l+1)\le n<g(k+2l+2)$ with $l\ge \mathcal{T}\{H_4,H_5,H_6\}$.
\\[10pt]
{\noindent\bf Case 4} $g(k+2l+2)\le n<f(k+2l+2)$ with $l\ge 1$. We have
\begin{align*}
\gamma_{Y,A,B}^k (n)&\ge \sum_{j=0}^{l}F(n-f(k+2j))-\sum_{j=0}^{l}F(n-g(k+2j+1))\\[5pt]
&-\sum_{j=0}^{l}F(n-f(k+2j+1))+\sum_{j=0}^{l}F(n-g(k+2j+2))-(4l+4)e\\[5pt]
&=P_6(n).
\end{align*}
Via Maple, we find that $P_6(n)$ is a quadratic polynomial in variable $n$ with
$LC_n(P_6(n))=-6a(l+1)(A-B)<0$. It follows that
\begin{align*}
P_6(n)\ge \min\{P_6(g(k+2l+2)),P_6(f(k+2l+2))\}.
\end{align*}
Let $H_7=P_6(g(k+2l+2))$ and $H_8=P_6(f(k+2l+2))$. Since $LC_l(H_7)=LC_l(H_8)=48A^2ak(A-B)>0$, by Lemma \ref{lem-3} we have $\gamma_{Y,A,B}^k (n)\ge P_6(n)\ge \min\{H_7,H_8\}\ge 0$ for
$l\ge \mathcal{T}_l\{H_7,H_8\}$.

Then $\gamma_{Y,A,B}^k (n)\ge 0$ for $g(k+2l+2)\le n<f(k+2l+2)$ with $l\ge \mathcal{T}_l\{H_7,H_8\}$.

Finally, combining the above four cases, we have $\gamma_{Y,A,B}^k (n)\ge 0$ for $n\ge f(k+2l)$ with $l\ge \lceil\mathcal{T}_l\{H_1,H_2,\cdots,H_8\}\rceil$.
The eight polynomials $H_1,H_2,\cdots,H_8$ in variable $l$ can be computed by Maple, which are listed in
Appendix A. This completes the proof of Theorem \ref{t-1}.

\section{Proof of Theorem \ref{t-2}}
Let $a,b,c,d$ be the quantities defined in Lemma \ref{lem-2} and $e=D_{z_1,z_2,z_3,z_4,z_5}$.
Let $F(x)=ax^4+bx^3+cx^2+dx$, $f(j)=Aj^2+Bj$ and $g(j)=Aj^2-Bj$.

Next, we shall distinguish four cases to prove Theorem \ref{t-2}.

{\noindent\bf Case 1} $f(k+2l)\le n<g(k+2l+1)$ with $l\ge 1$. We have
\begin{align*}
\gamma_{Z,A,B}^k (n)&\ge \sum_{j=0}^{l-1}F(n-f(k+2j))-\sum_{j=0}^{l-1}F(n-g(k+2j+1))\\[5pt]
&-\sum_{j=0}^{l-1}F(n-f(k+2j+1))+\sum_{j=0}^{l-1}F(n-g(k+2j+2))\\[5pt]
&+F(n-f(k+2l))-(4l+1)e\\[5pt]
&=Q_1(n).
\end{align*}
Via Maple, we find that $Q_1(n)$ is a polynomial of degree $4$ in variable $n$ with
$LC_n(Q_1(n))=a>0$ and
\begin{align*}
Q_2(n)=\frac{\mathrm{d}^2}{\mathrm{d}n^2}Q_1(n)=12an^2+Un+V,
\end{align*}
where $12a>0$, $U$ and $V$ are independent of $n$. It follows that
\begin{align*}
Q_2(n)\ge Q_2\left(-\frac{U}{24a}\right)=G_1.
\end{align*}
Note that $LC_l(G_1)=48a(A-B)(2Ak-A-B)>0$. By Lemma \ref{lem-3}, we conclude that for $l\ge \mathcal{T}_l\{G_1\}$, we have $Q_2(n)\ge G_1\ge 0$, and so $Q_3(n)=\frac{\mathrm{d}}{\mathrm{d}n}Q_1(n)$ is a weakly increasing function on the interval $[f(k+2l),g(k+2l+1))$.
It follows that for $l\ge \mathcal{T}_l\{G_1\}$,
\begin{align*}
Q_3(n)\ge Q_3(f(k+2l))=G_2.
\end{align*}
Note that $LC_l(G_2)=192A^2ak(A-B)>0$. By Lemma \ref{lem-3}, for $l\ge \mathcal{T}_l\{G_1,G_2\}$ we have $Q_3(n)\ge G_2\ge 0$, and so $Q_1(n)$ is a weakly increasing function on the interval $[f(k+2l),g(k+2l+1))$. It follows that for $l\ge \mathcal{T}_l\{G_1,G_2\}$,
\begin{align*}
\gamma_{Z,A,B}^k (n)\ge Q_1(n)\ge Q_1(f(k+2l))=G_3.
\end{align*}
Since $LC_l(G_3)=256 A^3 a k (A-B)>0$, by Lemma \ref{lem-3} we have $\gamma_{Z,A,B}^k (n)\ge Q_1(n)\ge G_3\ge 0$ for $l\ge \mathcal{T}_l\{G_1,G_2,G_3\}$.

Then $\gamma_{Z,A,B}^k (n)\ge 0$ for $f(k+2l)\le n<g(k+2l+1)$ with $l\ge \mathcal{T}_l\{G_1,G_2,G_3\}$.
\\[10pt]
{\noindent\bf Case 2} $g(k+2l+1)\le n<f(k+2l+1)$ with $l\ge 1$. We have
\begin{align*}
\gamma_{Z,A,B}^k (n)&\ge \sum_{j=0}^{l-1}F(n-f(k+2j))-\sum_{j=0}^{l-1}F(n-g(k+2j+1))\\[5pt]
&-\sum_{j=0}^{l-1}F(n-f(k+2j+1))+\sum_{j=0}^{l-1}F(n-g(k+2j+2))\\[5pt]
&+F(n-f(k+2l))-F(n-g(k+2l+1))-(4l+2)e\\[5pt]
&=Q_4(n).
\end{align*}
Via Maple, we find that $Q_4(n)$ is a cubic polynomial in variable $n$ with
$LC_n(Q_4(n))=4a(2k+2l+1)(A-B)>0$ and
\begin{align*}
Q_5(n)=\frac{\mathrm{d}}{\mathrm{d}n}Q_4(n)=12a(2k+2l+1)(A-B)n^2+Un+V,
\end{align*}
where $12a(2k+2l+1)(A-B)>0$, $U$ and $V$ are independent of $n$. It follows that
\begin{align*}
Q_5(n)\ge Q_5\left(-\frac{U}{24a(2k+2l+1)(A-B)}\right)=G_4.
\end{align*}
Note that $LC(G_4)=192A^2ak(A-B)>0$.
For $l\ge \mathcal{T}_l\{G_4\}$, we have $Q_5(n)\ge G_4\ge 0$, and so $Q_4(n)$ is a weakly increasing function on the interval $[g(k+2l+1),f(k+2l+1))$.
It follows that for $l\ge \mathcal{T}_l\{G_4\}$,
\begin{align*}
\gamma_{Z,A,B}^k (n)\ge Q_4(n)\ge Q_4(g(k+2l+1))=G_5.
\end{align*}
Since $LC_l(G_5)=256 A^3 a k (A-B)>0$, we have $\gamma_{Z,A,B}^k (n)\ge G_5\ge 0$ for $l\ge \mathcal{T}_l\{G_4,G_5\}$.

Then $\gamma_{Z,A,B}^k (n)\ge 0$ for $g(k+2l+1)\le n<f(k+2l+1)$ with $l\ge \mathcal{T}_l\{G_4,G_5\}$.
\\[10pt]
{\noindent\bf Case 3} $f(k+2l+1)\le n<g(k+2l+2)$ with $l\ge 1$. We have
\begin{align*}
\gamma_{Z,A,B}^k (n)&\ge \sum_{j=0}^{l-1}F(n-f(k+2j))-\sum_{j=0}^{l-1}F(n-g(k+2j+1))\\[5pt]
&-\sum_{j=0}^{l-1}F(n-f(k+2j+1))+\sum_{j=0}^{l-1}F(n-g(k+2j+2))\\[5pt]
&+F(n-f(k+2l))-F(n-g(k+2l+1))-F(n-f(k+2l+1))-(4l+3)e\\[5pt]
&=Q_6(n).
\end{align*}
Via Maple, we find that $Q_6(n)$ is a polynomial of degree $4$ in variable $n$ with
$LC_n(Q_6(n))=-a<0$ and
\begin{align*}
Q_7(n)=\frac{\mathrm{d}^2}{\mathrm{d}n^2}Q_6(n)=-12an^2+Un+V,
\end{align*}
where $-12a<0$, $U$ and $V$ are independent of $n$. It follows that
\begin{align*}
Q_7(n)\ge \min\{Q_7(f(k+2l+1)),Q_7(g(k+2l+2))\}.
\end{align*}
Let $G_6=Q_7(f(k+2l+1))$ and $G_7=Q_7(g(k+2l+2))$. Since $LC_l(G_6)=LC_l(G_7)=96 a (A-B) (A k+B)>0$,
for $l\ge \mathcal{T}_l\{G_6,G_7\}$ we have $Q_7(n)\ge \min\{G_6,G_7\}\ge 0$, and so $Q_8(n)=\frac{\mathrm{d}}{\mathrm{d}n}Q_6(n)$ is a weakly increasing function on the interval $[f(k+2l+1),g(k+2l+2))$.
It follows that for $l\ge \mathcal{T}_l\{G_6,G_7\}$,
\begin{align*}
Q_8(n)\ge Q_8(f(k+2l+1))=G_8.
\end{align*}
Since $LC_l(G_8)=192 A^2 a k (A-B)>0$, for $l\ge \mathcal{T}_l\{G_6,G_7,G_8\}$ we have $Q_8(n)\ge G_8\ge 0$, and so $Q_6(n)$ is a weakly increasing function on the interval $[f(k+2l+1),g(k+2l+2))$. It follows that for $l\ge \mathcal{T}_l\{G_6,G_7,G_8\}$,
\begin{align*}
\gamma_{Z,A,B}^k (n)\ge Q_6(n)\ge Q_6(f(k+2l+1))=G_9.
\end{align*}
Since $LC(G_9)=256 A^3 a k (A-B)>0$, we have $\gamma_{Z,A,B}^k (n)\ge G_9\ge 0$ for $l\ge \mathcal{T}_l\{G_6,G_7,G_8,G_9\}$.

Then we have $\gamma_{Z,A,B}^k (n)\ge 0$ for $f(k+2l+1)\le n<g(k+2l+2)$ with $l\ge \mathcal{T}_l\{G_6,G_7,G_8,G_9\}$.
\\[10pt]
{\noindent\bf Case 4} $g(k+2l+2)\le n<f(k+2l+2)$ with $l\ge 1$. We have
\begin{align*}
\gamma_{Z,A,B}^k (n)&\ge \sum_{j=0}^{l}F(n-f(k+2j))-\sum_{j=0}^{l}F(n-g(k+2j+1))\\[5pt]
&-\sum_{j=0}^{l}F(n-f(k+2j+1))+\sum_{j=0}^{l}F(n-g(k+2j+2))-(4l+4)e\\[5pt]
&=Q_9(n).
\end{align*}
Via Maple, we find that $Q_9(n)$ is a cubic polynomial in variable $n$ with
$LC_n(Q_9(n))=-8a(l+1)(A-B)<0$ and
\begin{align*}
Q_{10}(n)=\frac{\mathrm{d}}{\mathrm{d}n}Q_9(n)=-24a(l+1)(A-B)n^2+Un+V,
\end{align*}
where $-24a(l+1)(A-B)<0$, $U$ and $V$ are independent of $n$. It follows that
\begin{align*}
Q_{10}(n)\ge \min\{Q_{10}(g(k+2l+2)),Q_{10}(f(k+2l+2))\}.
\end{align*}
Let $G_{10}=Q_{10}(g(k+2l+2))$ and $G_{11}=Q_{10}(f(k+2l+2))$.
Since $LC_l(G_{10})=LC_l(G_{11})=192 A^2 a k (A-B)>0$, for $l\ge \mathcal{T}_l\{G_{10},G_{11}\}$
we have $Q_{10}(n)\ge \min\{G_{10},G_{11}\}\ge 0$, and so $Q_9(n)$ is a weakly increasing function on the interval $[g(k+2l+2),f(k+2l+2))$. It follows that for $l\ge \mathcal{T}_l\{G_{10},G_{11}\}$,
\begin{align*}
\gamma_{Z,A,B}^k (n)\ge Q_9(n)\ge Q_9(g(k+2l+2))=G_{12}.
\end{align*}
Since $LC_l(G_{12})=256 A^3 a k (A-B)>0$, we have $\gamma_{Z,A,B}^k (n)\ge G_{12}\ge 0$ for $l\ge \mathcal{T}_l\{G_{10},G_{11},G_{12}\}$.

Then $\gamma_{Z,A,B}^k (n)\ge 0$ for $g(k+2l+2)\le n<f(k+2l+2)$ with $l\ge \mathcal{T}_l\{G_{10},G_{11},G_{12}\}$.

Finally, combining the above four cases, we have $\gamma_{Z,A,B}^k (n)\ge 0$ for $n\ge f(k+2l)$ with $l\ge \lceil\mathcal{T}_l\{G_1,G_2,\cdots,G_{12}\}\rceil$.
The twelve polynomials $G_1,G_2,\cdots,G_{12}$ in variable $l$ can be computed by Maple, which are listed in Appendix B. This completes the proof of Theorem \ref{t-2}.

\vskip 5mm \noindent{\bf Acknowledgments.}
The first author was supported by the National Natural Science Foundation of China (grant 12171370).

\section*{Appendix A}
\begin{align*}
&H_1=12a(A-B)(2Ak-A-B)l^2+6a(2k+1)(A-B)(2Ak-A-B)l+\frac{3ac-b^2}{3a},\\[5pt]
&H_2=48A^2ak(A-B)l^4+16a(A-B)(6A^2k^2+3ABk+A^2+AB-2B^2)l^3\\[5pt]
&+4 (A-B) (12 A^2 a k^3+18 A B a k^2+9 A^2 a k+9 A B a k-9 B^2 a k-3 A B a+2 A b k-3 B^2 a\\[5pt]
&-2 B b) l^2+(24 A^2 B a k^3-24 A B^2 a k^3+24 A^3 a k^2-36 A B^2 a k^2+12 B^3 a k^2-6 A^2 B a k\\[5pt]
&+8 A^2 b k^2-8 A B b k^2+6 B^3 a k-6 A^3 a+6 A B^2 a-4 A B b k+4 B^2 b k-2 A^2 b+2 B^2 b-2 A c\\[5pt]
&+2 B c-4 e)l-e,\\[5pt]
&H_3=48 A^2 a k (A-B) l^4+8 a (A-B) (9 A^2 k^2+12 A^2 k-2 A^2-2 B A+B^2) l^3\\[5pt]
&+12 a (A-B) (2 A^2 k^3+9 A^2 k^2+4 A^2 k-2 A B k+B^2 k-2 A^2-2 B A+B^2) l^2\\[5pt]
&+((144 A^3 a^2 k^3-144 A^2 B a^2 k^3+252 A^3 a^2 k^2-324 A^2 B a^2 k^2+108 A B^2 a^2 k^2-36 B^3 a^2 k^2\\[5pt]
&-144 A^2 B a^2 k+216 A B^2 a^2 k-72 B^3 a^2 k-45 A^3 a^2-9 A^2 B a^2+81 A B^2 a^2-27 B^3 a^2\\[5pt]
&+12 A c a-4 A b^2-12 B c a+4 B b^2-24 e a) l)/(6 a)
+(24A^3a^2k^3-72A^2Ba^2k^3\\[5pt]
&+72AB^2a^2k^3-24B^3a^2k^3+36A^3a^2k^2-108A^2Ba^2k^2+108AB^2a^2k^2-36B^3a^2k^2\\[5pt]
&+18A^3a^2k-54A^2Ba^2k+54AB^2a^2k-18B^3a^2k+3A^3a^2-9A^2Ba^2+9AB^2a^2-3B^3a^2\\[5pt]
&+24Aack-8Ab^2k-24Back+8Bb^2k+12Aac-4Ab^2-12Bac+4Bb^2-24ae)/(12a),\\[5pt]
&H_4=24a(A-B)(Ak+B)l^2+2(A-B)(12Aak^2+12Aak+18Bak+3Aa+15Ba+2b)l\\[5pt]
&+12A^2ak^2-12B^2ak^2+12A^2ak+12ABak-24B^2ak+3aA^2+6ABa+4Abk-9B^2a\\[5pt]
&-4Bbk+2Ab-2Bb-c,\\[5pt]
&H_5=24a(A-B)(Ak+B)l^2+2(A-B)(12Aak^2+24Aak+6Bak-3Aa+21Ba-2b)l\\[5pt]
&+24A^2ak^2-24ABak^2+24A^2ak-12ABak-12B^2ak-6A^2a+24ABa-18B^2a\\[5pt]
&-4Ab+4Bb-c,\\[5pt]
&H_6=48A^2ak(A-B)l^4+16 a (A-B) (6 A^2 k^2+6 A^2 k+3 A B k-A^2-A B+2 B^2) l^3\\[5pt]
&+4(A-B)(12A^2ak^3+36A^2ak^2+18ABak^2+15A^2ak+15ABak+15B^2ak-6A^2a\\[5pt]
&-3ABa+2Abk+15B^2a+2Bb)l^2+(48A^3ak^3-24A^2Bak^3-24AB^2ak^3+72A^3ak^2\\[5pt]
&-36AB^2ak^2-36B^3ak^2+12A^3ak+30A^2Bak+8A^2bk^2+36AB^2ak-8ABbk^2\\[5pt]
&-78B^3ak-6A^3a+12A^2Ba+8A^2bk+30AB^2a+4ABbk-36B^3a-12B^2bk+2A^2b\\[5pt]
&+8ABb-10B^2b+2Ac-2Bc-4e)l+8A^3ak^3-8B^3ak^3+12A^3ak^2+12A^2Bak^2\\[5pt]
&-24B^3ak^2+6A^3ak+12A^2Bak+4A^2bk^2+6AB^2ak-24B^3ak-4B^2bk^2+A^3a\\[5pt]
&+3A^2Ba+4A^2bk+3AB^2a+4ABbk-7B^3a-8B^2bk+A^2b+2ABb+2Ack-3B^2b\\[5pt]
&-2Bck+Ac-Bc-3e,
\end{align*}

\begin{align*}
&H_7=48A^2ak(A-B)l^4+16a(A-B)(6A^2k^2+12A^2k-3ABk+A^2+AB-2B^2)l^3\\[5pt]
&+4(A-B)(12A^2ak^3+72A^2ak^2-18ABak^2+81A^2ak-27ABak-9B^2ak\\[5pt]
&+12A^2a+15ABa+2Abk-21B^2a+2Bb)l^2+(96A^3ak^3-120A^2Bak^3\\[5pt]
&+24AB^2ak^3+312A^3ak^2-432A^2Bak^2+108AB^2ak^2+12B^3ak^2+264A^3ak\\[5pt]
&-330A^2Bak+8A^2bk^2-8ABbk^2+66B^3ak+42A^3a+24A^2Ba+16A^2bk\\[5pt]
&-138AB^2a-12ABbk+72B^3a-4B^2bk-2A^2b+16ABb-14B^2b-2Ac+2Bc\\[5pt]
&-4e)l+48A^3ak^3-72A^2Bak^3+24AB^2ak^3+120A^3ak^2-168A^2Bak^2+36AB^2ak^2\\[5pt]
&+12B^3ak^2+84A^3ak-90A^2Bak+8A^2bk^2-24AB^2ak-8ABbk^2+30B^3ak+10A^3a\\[5pt]
&+12A^2Ba+8A^2bk-42AB^2a-4ABbk+20B^3a-4B^2bk-2A^2b+8ABb-6B^2b\\[5pt]
&-2Ac+2Bc-4e,\\[5pt]
&H_8=48A^2ak(A-B)l^4+16a(A-B)(6A^2k^2+12A^2k+3ABk+A^2+AB-2B^2)l^3\\[5pt]
&+4(A-B)(12A^2ak^3+72A^2ak^2+18ABak^2+81A^2ak+45ABak-9B^2ak+12A^2a\\[5pt]
&+9ABa+2Abk-27B^2a-2Bb)l^2+(96A^3ak^3-72A^2Bak^3-24AB^2ak^3+312A^3ak^2\\[5pt]
&-144A^2Bak^2-180AB^2ak^2+12B^3ak^2+264A^3ak-54A^2Bak+8A^2bk^2-288AB^2ak\\[5pt]
&-8ABbk^2+78B^3ak+42A^3a-24A^2Ba+16A^2bk-138AB^2a-20ABbk+120B^3a\\[5pt]
&+4B^2bk-2A^2b-16ABb+18B^2b-2Ac+2Bc-4e)l+48A^3ak^3-24A^2Bak^3\\[5pt]
&-24AB^2ak^3+120A^3ak^2-24A^2Bak^2-108AB^2ak^2+12B^3ak^2+84A^3ak-6A^2Bak\\[5pt]
&+8A^2bk^2-120AB^2ak-8ABbk^2+42B^3ak+10A^3a-12A^2Ba+8A^2bk-42AB^2a\\[5pt]
&-12ABbk+44B^3a+4B^2bk-2A^2b-8ABb+10B^2b-2Ac+2Bc-4e,
\end{align*}
where $a,b,c$ and $e$ are given by
\begin{align*}
a&=\frac{1}{6y_1y_2y_3y_4},\\[5pt]
b&=\frac{y_1+y_2+y_3+y_4}{4y_1y_2y_3y_4},\\[5pt]
c&=\frac{y_1^2+y_2^2+y_3^2+y_4^2+3(y_1y_2+y_1y_3+y_1y_4+y_2y_3+y_2y_4+y_3y_4)}{12y_1y_2y_3y_4},\\[5pt]
e&=D_{y_1,y_2,y_3,y_4}.
\end{align*}

\section*{Appendix B}
\begin{align*}
&G_1=48a(A-B)(2Ak-A-B)l^2+24a(2k+1)(A-B)(2Ak-A-B)l+\frac{8ac-3b^2}{4a},\\[5pt]
&G_2=192A^2ak(A-B)l^4+64a(A-B)(6A^2k^2+3ABk+A^2+AB-2B^2)l^3\\[5pt]
&+24 (A-B) (8 A^2 a k^3+12 A B a k^2+6 A^2 a k+6 A B a k-6 B^2 a k-2 A B a+A b k-2 B^2 a\\[5pt]
&-B b) l^2+2 (A-B) (48 A B a k^3+48 A^2 a k^2+48 A B a k^2-24 B^2 a k^2-12 A B a k+12 A b k^2\\[5pt]
&-12 B^2 a k-12 A^2 a-12 A B a-6 B b k-3 A b-3 B b-2 c) l+d,\\[5pt]
&G_3=256 A^3 a k (A-B) l^6+384 A^2 a k (A-B) (2 A k+B) l^5+16 (A-B) (48 A^3 a k^3\\[5pt]
&+60 A^2 B a k^2+6 A^3 a k+6 A^2 B a k+12 A B^2 a k+8 A^2 B a+3 A^2 b k+8 A B^2 a-8 B^3 a) l^4\\[5pt]
&+16 (A-B) (16 A^3 a k^4+48 A^2 B a k^3+12 A^3 a k^2+12 A^2 B a k^2+24 A B^2 a k^2+22 A^2 B a k\\[5pt]
&+6 A^2 b k^2+22 A B^2 a k-14 B^3 a k+8 A^3 a+8 A^2 B a-4 A B^2 a+3 A B b k-4 B^3 a+A^2 b\\[5pt]
&+A B b-2 B^2 b) l^3+4 (A-B) (48 A^2 B a k^4+16 A^3 a k^3+16 A^2 B a k^3+64 A B^2 a k^3\\[5pt]
&+84 A^2 B a k^2+12 A^2 b k^3+84 A B^2 a k^2-36 B^3 a k^2+60 A^3 a k+60 A^2 B a k-18 A B^2 a k\\[5pt]
&+18 A B b k^2-18 B^3 a k-12 A^2 B a+9 A^2 b k-12 A B^2 a+9 A B b k-9 B^2 b k-3 A B b\\[5pt]
&+2 A c k-3 B^2 b-2 B c) l^2+(-32 A^4 a k^4+96 A^2 B^2 a k^4-64 A B^3 a k^4+96 A^3 B a k^3\\[5pt]
&-128 A B^3 a k^3+32 B^4 a k^3+144 A^4 a k^2-168 A^2 B^2 a k^2+24 A^2 B b k^3-24 A B^2 b k^3\\[5pt]
&+24 B^4 a k^2-24 A^3 B a k+24 A^3 b k^2+24 A B^3 a k-36 A B^2 b k^2+12 B^3 b k^2-34 A^4 a\\[5pt]
&+36 A^2 B^2 a-6 A^2 B b k+8 A^2 c k^2-8 A B c k^2-2 B^4 a+6 B^3 b k-6 A^3 b+6 A B^2 b-4 A B c k\\[5pt]
&+4 B^2 c k-2 A^2 c+2 B^2 c-2 A d+2 B d-4 e) l-e,\\[5pt]
&G_4=192A^2ak(A-B)l^4+32 a (A-B) (9 A^2 k^2+12 A^2 k-2 A^2-2 A B+B^2) l^3\\[5pt]
&+48 a (A-B) (2 A^2 k^3+9 A^2 k^2+4 A^2 k-2 A B k+B^2 k-2 A^2-2 A B+B^2) l^2\\[5pt]
&+((A-B)(192 A^2 a^2 k^3+336 A^2 a^2 k^2-96 A B a^2 k^2+48 B^2 a^2 k^2-192 A B a^2 k\\[5pt]
&+96 B^2 a^2 k-60 A^2 a^2-72 A B a^2+36 B^2 a^2+8 a c-3 b^2) l)/(2 a)\\[5pt]
&+((2 k+1) (A-B) (16 A^2 a^2 k^2-32 A B a^2 k^2+16 B^2 a^2 k^2+16 A^2 a^2 k-32 A B a^2 k\\[5pt]
&+16 B^2 a^2 k+4 A^2 a^2-8 A B a^2+4 B^2 a^2+8 a c-3 b^2))/(4 a),\\[5pt]
\end{align*}

\begin{align*}
&G_5=256 A^3 a k (A-B) l^6+384 A^2 a k (A-B) (2 k A+2 A-B) l^5\\[5pt]
&+16 (A-B) (48 A^3 a k^3+120 A^3 a k^2-60 A^2 B a k^2+66 A^3 a k-54 A^2 B a k+12 A B^2 a k\\[5pt]
&+8 A^2 B a+3 A^2 b k+8 A B^2 a-8 B^3 a) l^4+16 (A-B) (16 A^3 a k^4+96 A^3 a k^3-48 A^2 B a k^3\\[5pt]
&+132 A^3 a k^2-108 A^2 B a k^2+24 A B^2 a k^2+52 A^3 a k-38 A^2 B a k+6 A^2 b k^2+34 A B^2 a k\\[5pt]
&-18 B^3 a k-8 A^3 a+8 A^2 B a+6 A^2 b k+20 A B^2 a-3 A B b k-12 B^3 a-A^2 b-A B b\\[5pt]
&+2 B^2 b) l^3+4 (A-B) (96 A^3 a k^4-48 A^2 B a k^4+304 A^3 a k^3-272 A^2 B a k^3+64 A B^2 a k^3\\[5pt]
&+312 A^3 a k^2-276 A^2 B a k^2+12 A^2 b k^3+156 A B^2 a k^2-60 B^3 a k^2+60 A^3 a k-60 A^2 B a k\\[5pt]
&+36 A^2 b k^2+162 A B^2 a k-18 A B b k^2-78 B^3 a k-48 A^3 a-12 A^2 B a+15 A^2 b k+60 A B^2 a\\[5pt]
&-21 A B b k-24 B^3 a+15 B^2 b k-6 A^2 b-9 A B b+2 A c k+9 B^2 b-2 B c) l^2+(160 A^4 a k^4\\[5pt]
&-384 A^3 B a k^4+288 A^2 B^2 a k^4-64 A B^3 a k^4+448 A^4 a k^3-928 A^3 B a k^3+768 A^2 B^2 a k^3\\[5pt]
&-384 A B^3 a k^3+96 B^4 a k^3+336 A^4 a k^2-672 A^3 B a k^2+48 A^3 b k^3+744 A^2 B^2 a k^2-72 A^2 B b k^3\\[5pt]
&-576 A B^3 a k^2+24 A B^2 b k^3+168 B^4 a k^2-48 A^4 a k-120 A^3 B a k+72 A^3 b k^2+432 A^2 B^2 a k\\[5pt]
&-144 A^2 B b k^2-360 A B^3 a k+108 A B^2 b k^2+96 B^4 a k-36 B^3 b k^2-62 A^4 a+16 A^3 B a+12 A^3 b k\\[5pt]
&+108 A^2 B^2 a-78 A^2 B b k+8 A^2 c k^2-80 A B^3 a+108 A B^2 b k-8 A B c k^2+18 B^4 a-42 B^3 b k\\[5pt]
&-6 A^3 b-12 A^2 B b+8 A^2 c k+30 A B^2 b-20 A B c k-12 B^3 b+12 B^2 c k+2 A^2 c-8 A B c+6 B^2 c\\[5pt]
&+2 A d-2 B d-4 e) l-128 A B^3 a k^3-96 A^3 B a k^2+144 A^2 B^2 a k^2-24 A^2 B b k^3-96 A B^3 a k^2\\[5pt]
&-32 A^3 B a k+48 A^2 B^2 a k-36 A^2 B b k^2-32 A B^3 a k+36 A B^2 b k^2-18 A^2 B b k+18 A B^2 b k\\[5pt]
&-8 A B c k^2-8 A B c k-64 A^3 B a k^4+96 A^2 B^2 a k^4-64 A B^3 a k^4-128 A^3 B a k^3+192 A^2 B^2 a k^3\\[5pt]
&-2 e,\\[5pt]
&G_6=96 a (A-B) (A k+B) l^2+12 (A-B) (8 A a k^2+8 A a k+12 B a k+2 A a+10 B a+b) l\\[5pt]
&+48 A^2 a k^2-48 B^2 a k^2+48 A^2 a k+48 A B a k-96 B^2 a k+12 A^2 a+24 A B a+12 A b k\\[5pt]
&-36 B^2 a-12 B b k+6 A b-6 B b-2 c,\\[5pt]
&G_7=96 a (A-B) (A k+B) l^2+12 (A-B) (8 A a k^2+16 A a k+4 B a k-2 A a+14 B a-b) l\\[5pt]
&+96 A^2 a k^2-96 A B a k^2+96 A^2 a k-48 A B a k-48 B^2 a k-24 A^2 a+96 A B a-72 B^2 a\\[5pt]
&-12 A b+12 B b-2 c,
\end{align*}

\begin{align*}
&G_8=192 A^2 a k (A-B) l^4+64 a (A-B) (6 A^2 k^2+6 A^2 k+3 A B k-A^2-B A+2 B^2) l^3\\[5pt]
&+24 (A-B) (8 A^2 a k^3+24 A^2 a k^2+12 A B a k^2+10 A^2 a k+10 A B a k+10 B^2 a k-4 A^2 a\\[5pt]
&-2 A B a+A b k+10 B^2 a+B b) l^2+2 (A-B) (96 A^2 a k^3+48 A B a k^3+144 A^2 a k^2\\[5pt]
&+144 A B a k^2+72 B^2 a k^2+24 A^2 a k+84 A B a k+12 A b k^2+156 B^2 a k-12 A^2 a+12 A B a\\[5pt]
&+12 A b k+72 B^2 a+18 B b k+3 A b+15 B b+2 c) l+32 A^3 a k^3-32 B^3 a k^3+48 A^3 a k^2\\[5pt]
&+48 A^2 B a k^2-96 B^3 a k^2+24 A^3 a k+48 A^2 B a k+12 A^2 b k^2+24 A B^2 a k-96 B^3 a k\\[5pt]
&-12 B^2 b k^2+4 A^3 a+12 A^2 B a+12 A^2 b k+12 A B^2 a+12 A B b k-28 B^3 a-24 B^2 b k\\[5pt]
&+3 A^2 b+6 A B b+4 A c k-9 B^2 b-4 B c k+2 A c-2 B c-d,\\[5pt]
&G_9=256 A^3 a k (A-B) l^6+384 A^2 a k (A-B) (2 A k+2 A+B) l^5\\[5pt]
&+16 (A-B) (48 A^3 a k^3+120 A^3 a k^2+60 A^2 B a k^2+66 A^3 a k+66 A^2 B a k+12 A B^2 a k\\[5pt]
&-8 A^2 B a+3 A^2 b k-8 A B^2 a+8 B^3 a) l^4+16 (A-B) (16 A^3 a k^4+96 A^3 a k^3+48 A^2 B a k^3\\[5pt]
&+132 A^3 a k^2+132 A^2 B a k^2+24 A B^2 a k^2+52 A^3 a k+62 A^2 B a k+6 A^2 b k^2+14 A B^2 a k\\[5pt]
&+18 B^3 a k-8 A^3 a-24 A^2 B a+6 A^2 b k-12 A B^2 a+3 A B b k+20 B^3 a-A^2 b-A B b\\[5pt]
&+2 B^2 b) l^3+4 (A-B) (96 A^3 a k^4+48 A^2 B a k^4+304 A^3 a k^3+304 A^2 B a k^3+64 A B^2 a k^3\\[5pt]
&+312 A^3 a k^2+420 A^2 B a k^2+12 A^2 b k^3+132 A B^2 a k^2+60 B^3 a k^2+60 A^3 a k+60 A^2 B a k\\[5pt]
&+36 A^2 b k^2+42 A B^2 a k+18 A B b k^2+138 B^3 a k-48 A^3 a-84 A^2 B a+15 A^2 b k-12 A B^2 a\\[5pt]
&+15 A B b k+72 B^3 a+15 B^2 b k-6 A^2 b-3 A B b+2 A c k+15 B^2 b+2 B c) l^2+(160 A^4 a k^4\\[5pt]
&-96 A^2 B^2 a k^4-64 A B^3 a k^4+448 A^4 a k^3+160 A^3 B a k^3-384 A^2 B^2 a k^3-128 A B^3 a k^3\\[5pt]
&-96 B^4 a k^3+336 A^4 a k^2+192 A^3 B a k^2+48 A^3 b k^3-216 A^2 B^2 a k^2-24 A^2 B b k^3-24 A B^2 b k^3\\[5pt]
&-312 B^4 a k^2-48 A^4 a k+24 A^3 B a k+72 A^3 b k^2+192 A^2 B^2 a k+168 A B^3 a k-36 A B^2 b k^2\\[5pt]
&-336 B^4 a k-36 B^3 b k^2-62 A^4 a-16 A^3 B a+12 A^3 b k+108 A^2 B^2 a+30 A^2 B b k+8 A^2 c k^2\\[5pt]
&+80 A B^3 a+36 A B^2 b k-8 A B c k^2-110 B^4 a-78 B^3 b k-6 A^3 b+12 A^2 B b+8 A^2 c k\\[5pt]
&+30 A B^2 b+4 A B c k-36 B^3 b-12 B^2 c k+2 A^2 c+8 A B c-10 B^2 c+2 A d-2 B d-4 e) l\\[5pt]
&+16 A^4 a k^4-16 B^4 a k^4+32 A^4 a k^3+32 A^3 B a k^3-64 B^4 a k^3+24 A^4 a k^2+48 A^3 B a k^2\\[5pt]
&+8 A^3 b k^3+24 A^2 B^2 a k^2-96 B^4 a k^2-8 B^3 b k^3+8 A^4 a k+24 A^3 B a k+12 A^3 b k^2+24 A^2 B^2 a k\\[5pt]
&+12 A^2 B b k^2+8 A B^3 a k-64 B^4 a k-24 B^3 b k^2+A^4 a+4 A^3 B a+6 A^3 b k+6 A^2 B^2 a+12 A^2 B b k\\[5pt]
&+4 A^2 c k^2+4 A B^3 a+6 A B^2 b k-15 B^4 a-24 B^3 b k-4 B^2 c k^2+A^3 b+3 A^2 B b+4 A^2 c k\\[5pt]
&+3 A B^2 b+4 A B c k-7 B^3 b-8 B^2 c k+A^2 c+2 A B c+2 A d k-3 B^2 c-2 B d k+A d-B d-3 e,
\end{align*}

\begin{align*}
&G_{10}=192 A^2 a k (A-B) l^4+64 a (A-B) (6 A^2 k^2+12 A^2 k-3 A B k+A^2+B A-2 B^2) l^3\\[5pt]
&+24 (A-B) (8 A^2 a k^3+48 A^2 a k^2-12 A B a k^2+54 A^2 a k-18 A B a k-6 B^2 a k\\[5pt]
&+8 A^2 a+10 A B a+A b k-14 B^2 a+B b) l^2+2 (A-B) (192 A^2 a k^3-48 A B a k^3+624 A^2 a k^2\\[5pt]
&-240 A B a k^2-24 B^2 a k^2+528 A^2 a k-132 A B a k+12 A b k^2-132 B^2 a k+84 A^2 a+132 A B a\\[5pt]
&+24 A b k-144 B^2 a+6 B b k-3 A b+21 B b-2 c) l+2 (A-B) (96 A^2 a k^3-48 A B a k^3\\[5pt]
&+240 A^2 a k^2-96 A B a k^2-24 B^2 a k^2+168 A^2 a k-12 A B a k+12 A b k^2-60 B^2 a k+20 A^2 a\\[5pt]
&+44 A B a+12 A b k-40 B^2 a+6 B b k-3 A b+9 B b-2 c),\\[5pt]
&G_{11}=192 A^2 a k (A-B) l^4+64 a (A-B) (6 A^2 k^2+12 A^2 k+3 A B k+A^2+B A-2 B^2) l^3\\[5pt]
&+24 (A-B) (8 A^2 a k^3+48 A^2 a k^2+12 A B a k^2+54 A^2 a k+30 A B a k-6 B^2 a k+8 A^2 a\\[5pt]
&+6 A B a+A b k-18 B^2 a-B b) l^2+2 (A-B) (192 A^2 a k^3+48 A B a k^3+624 A^2 a k^2\\[5pt]
&+336 A B a k^2-24 B^2 a k^2+528 A^2 a k+420 A B a k+12 A b k^2-156 B^2 a k+84 A^2 a+36 A B a\\[5pt]
&+24 A b k-240 B^2 a-6 B b k-3 A b-27 B b-2 c) l+2 (A-B) (96 A^2 a k^3+48 A B a k^3\\[5pt]
&+240 A^2 a k^2+192 A B a k^2-24 B^2 a k^2+168 A^2 a k+156 A B a k+12 A b k^2-84 B^2 a k\\[5pt]
&+20 A^2 a-4 A B a+12 A b k-88 B^2 a-6 B b k-3 A b-15 B b-2 c),
\end{align*}

\begin{align*}
&G_{12}=256 A^3 a k (A-B) l^6+384 A^2 a k (A-B) (2 A k+4 A-B) l^5\\[5pt]
&+16 (A-B) (48 A^3 a k^3+240 A^3 a k^2-60 A^2 B a k^2+246 A^3 a k-114 A^2 B a k+12 A B^2 a k\\[5pt]
&-8 A^2 B a+3 A^2 b k-8 A B^2 a+8 B^3 a) l^4+16 (A-B) (16 A^3 a k^4+192 A^3 a k^3-48 A^2 B a k^3\\[5pt]
&+492 A^3 a k^2-228 A^2 B a k^2+24 A B^2 a k^2+344 A^3 a k-238 A^2 B a k+6 A^2 b k^2+26 A B^2 a k\\[5pt]
&+14 B^3 a k+8 a A^3-24 A^2 B a+12 A^2 b k-36 A B^2 a-3 A B b k+28 B^3 a+A^2 b+A B b\\[5pt]
&-2 B^2 b) l^3+4 (A-B) (192 A^3 a k^4-48 A^2 B a k^4+1168 A^3 a k^3-560 A^2 B a k^3+64 A B^2 a k^3\\[5pt]
&+2064 A^3 a k^2-1380 A^2 B a k^2+12 A^2 b k^3+204 A B^2 a k^2+36 B^3 a k^2+1164 A^3 a k\\[5pt]
&-1020 A^2 B a k+72 A^2 b k^2+6 A B^2 a k-18 A B b k^2+150 B^3 a k+96 a A^3-84 A^2 B a\\[5pt]
&+81 A^2 b k-228 A B^2 a-27 A B b k+144 B^3 a-9 B^2 b k+12 A^2 b+15 A B b+2 A c k-21 B^2 b\\[5pt]
&+2 B c) l^2+(736 A^4 a k^4-1152 A^3 B a k^4+480 A^2 B^2 a k^4-64 A B^3 a k^4+3200 A^4 a k^3\\[5pt]
&-5472 A^3 B a k^3+2688 A^2 B^2 a k^3-384 A B^3 a k^3-32 B^4 a k^3+4560 A^4 a k^2-8352 A^3 B a k^2\\[5pt]
&+96 A^3 b k^3+4248 A^2 B^2 a k^2-120 A^2 B b k^3-192 A B^3 a k^2+24 A B^2 b k^3-264 B^4 a k^2\\[5pt]
&+2400 A^4 a k-4488 A^3 B a k+312 A^3 b k^2+1680 A^2 B^2 a k-432 A^2 B b k^2+936 A B^3 a k\\[5pt]
&+108 A B^2 b k^2-528 B^4 a k+12 B^3 b k^2+350 A^4 a-416 A^3 B a+264 A^3 b k-540 A^2 B^2 a\\[5pt]
&-330 A^2 B b k+8 A^2 c k^2+928 A B^3 a-8 A B c k^2-322 B^4 a+66 B^3 b k+42 A^3 b+24 A^2 B b\\[5pt]
&+16 A^2 c k-138 A B^2 b-12 A B c k+72 B^3 b-4 B^2 c k-2 A^2 c+16 A B c-14 B^2 c-2 A d\\[5pt]
&+2 B d-4 e) l+224 A^4 a k^4-448 A^3 B a k^4+288 A^2 B^2 a k^4-64 A B^3 a k^4+832 A^4 a k^3\\[5pt]
&-1632 A^3 B a k^3+960 A^2 B^2 a k^3-128 A B^3 a k^3-32 B^4 a k^3+1104 A^4 a k^2-2064 A^3 B a k^2\\[5pt]
&+48 A^3 b k^3+984 A^2 B^2 a k^2-72 A^2 B b k^3+96 A B^3 a k^2+24 A B^2 b k^3-120 B^4 a k^2+592 A^4 a k\\[5pt]
&-968 A^3 B a k+120 A^3 b k^2+168 A^2 B^2 a k-168 A^2 B b k^2+360 A B^3 a k+36 A B^2 b k^2-152 B^4 a k\\[5pt]
&+12 B^3 b k^2+94 A^4 a-80 A^3 B a+84 A^3 b k-156 A^2 B^2 a-90 A^2 B b k+8 A^2 c k^2+208 A B^3 a\\[5pt]
&-24 A B^2 b k-8 A B c k^2-66 B^4 a+30 B^3 b k+10 A^3 b+12 A^2 B b+8 A^2 c k-42 A B^2 b-4 A B c k\\[5pt]
&+20 B^3 b-4 B^2 c k-2 A^2 c+8 A B c-6 B^2 c-2 A d+2 B d-4 e,
\end{align*}
where $a,b,c,d$ and $e$ are given by
\begin{align*}
a&=\frac{1}{24z_1z_2z_3z_4z_5},\\[5pt]
b&=\frac{z_1+z_2+z_3+z_4+z_5}{12z_1z_2z_3z_4z_5},\\[5pt]
c&=\frac{z_1^2+z_2^2+z_3^2+z_4^2+z_5^2}{24z_1z_2z_3z_4z_5}\\[5pt]
&+\frac{z_1z_2+z_1z_3+z_1z_4+z_1z_5+z_2z_3+z_2z_4+z_2z_5+z_3z_4+z_3z_5+z_4z_5}{8z_1z_2z_3z_4z_5},\\[5pt]
d&=\frac{(z_1+z_2+z_3+z_4+z_5)}{24z_1z_2z_3z_4z_5}\\[5pt]
&\times (z_1z_2+z_1z_3+z_1z_4+z_1z_5+z_2z_3+z_2z_4+z_2z_5+z_3z_4+z_3z_5+z_4z_5),\\[5pt]
e&=D_{z_1,z_2,z_3,z_4,z_5}.
\end{align*}

\end{document}